\documentclass[11pt]{amsart}
\usepackage{amsmath,amsthm, amscd, amssymb, amsfonts}
\usepackage[all]{xy}
\usepackage{mathrsfs}
\usepackage{enumitem}

\usepackage[ansinew]{inputenc}
\usepackage{graphicx,fancyhdr}

\usepackage[dvips, dvipsnames, usenames]{color}

\newcommand{\comment}[1]{}

\numberwithin{equation}{section}
\newtheorem{theorem}{Theorem}[section]
\newtheorem{lemma}[theorem]{Lemma}

\newtheorem{prop}[theorem]{Proposition}

\theoremstyle{definition}

\newtheorem{question}[theorem]{Question}

\theoremstyle{remark}
\newtheorem{remark}[theorem]{Remark}

\def\pf{\begin{proof}}
\def\epf{\end{proof}}

\newcommand{\kk}{ \mathbf{k}}
\newcommand{\ku}{ \Bbbk}
\newcommand{\kut}{ \ku^{\times}}

\newcommand\G{\mathbb G}

\newcommand\ghost{\mathscr{G}}

\newcommand\I{\mathbb I}

\newcommand\N{\mathbb N}

\renewcommand{\_}[1]{_{\left( #1 \right)}}

\newcommand\cA{\mathcal{A}}
\newcommand\cB{\mathcal{B}}

\newcommand\D{\mathcal{D}}
\newcommand\E{\mathcal{E}}

\newcommand{\Ss}{{\mathcal S}}
\newcommand\T{\mathcal{T}}
\newcommand\cV{\mathcal{V}}

\newcommand\Ug{\mathfrak U}

\newcommand\ad{\operatorname{ad}}
\newcommand{\Alg}{\Hom_{\text{alg}}}
\newcommand\Aut{\operatorname{Aut}}
\newcommand{\AuH}{\Aut_{\text{Hopf}}}

\newcommand\car{\operatorname{char}}
\newcommand\Der{\operatorname{Der}}

\newcommand\id{\operatorname{id}}
\newcommand\gr{\operatorname{gr}}
\newcommand\GK{\operatorname{GKdim}}
\newcommand{\Hom}{\operatorname{Hom}}

\newcommand\soc{\operatorname{soc}}

\def\ydk{{}^{K}_{K}\mathcal{YD}}

\newcommand{\NA}{\mathcal{B}}
\newcommand{\toba}{\mathcal{B}}
\newcommand{\ot}{\otimes}

\newcommand{\ydG}{{}^{\ku \Gamma }_{\ku \Gamma }\mathcal{YD}}
\newcommand{\ydg}{{}^{\ku G}_{\ku G}\mathcal{YD}}

\newcounter{tabla}\stepcounter{tabla}

\begin{document}

\title[Liftings of Jordan and super Jordan planes]{Liftings of Jordan and super Jordan planes}

\author[Andruskiewitsch; Angiono; Heckenberger]
{Nicol\'as Andruskiewitsch, Iv\'an Angiono, Istv\'an Heckenberger}

\address{FaMAF-CIEM (CONICET), Universidad Nacional de C\'ordoba,
Medina A\-llen\-de s/n, Ciudad Universitaria, 5000 C\' ordoba, Rep\'
ublica Argentina.} \email{(andrus|angiono)@mate.uncor.edu}

\address{Philipps-Universität Marburg,
	Fachbereich Mathematik und Informatik,
	Hans-Meerwein-Straße,
	D-35032 Marburg, Germany.} \email{heckenberger@mathematik.uni-marburg.de}

\thanks{\noindent 2000 \emph{Mathematics Subject Classification.}
16T05. \newline The work of N. A. and I. A. was partially supported by CONICET,
FONCyT-ANPCyT, Secyt (UNC), the MathAmSud project GR2HOPF. The work of N. A., respectively I. A., was partially done during a visit to the University of Marburg, respectively the MPI (Bonn), supported by the Alexander von Humboldt Foundation. }

\begin{abstract}
We classify pointed Hopf algebras with finite Gelfand-Kirillov dimension whose infinitesimal braiding has dimension 2 but is not of diagonal type, or equivalently is a block. These Hopf algebras are new and turn out to be liftings of either a Jordan or a super Jordan plane over a nilpotent-by-finite group.
\end{abstract}

\maketitle

\setcounter{tocdepth}{1}

\section{Introduction}\label{section:introduction}
Let $\ku = \overline{\ku}$ be a field, $\car \ku = 0$. Let $H$ be a pointed Hopf algebra, $G = G(H)$, $\gr H$ the graded Hopf algebra associated to its coradical filtration, $R = \oplus_{n\ge 0} R^n$ the graded Hopf algebra in the category $\ydg$ of Yetter-Drinfeld modules such that $\gr H \simeq R \# \ku G$ and $V = R^1$ the infinitesimal braiding of $H$.
The classification of Hopf algebras with finite Gelfand-Kirillov dimension ($\GK$ for short) attracted considerable interest recently, see \cite{Ob}. Connected Hopf algebras with finite Gelfand-Kirillov dimension are quantum deformations of algebraic unipotent groups \cite{EG}. Also, there are several results in low $\GK$, see \cite{BZ, GZ, WZZ} and references therein.
Further, the  classification is known assuming that $H$ is a domain, $G$ is abelian and $V$ is of diagonal type \cite{AA,AS-crelle}. Here we contribute  to this question.

\medbreak
Let $\ell \in\N_{\ge 2}$ and  $\I_\ell=\{1,2,\dots,\ell\}$.
Let $\epsilon \in \kut$.
Let  $\cV(\epsilon, \ell)$ be the braided vector space with basis $(x_i)_{i\in
\I_\ell}$ and braiding $c\in \Aut (V\otimes V)$ such that
\begin{align}\label{eq:bloque2}
c(x_i \otimes x_1) &= \epsilon x_1 \otimes x_i,& c(x_i \otimes x_j) &= (\epsilon x_j + x_{j-1}) \otimes x_i,& i,j &\in \I_\ell.
\end{align}
We say that a braided vector space is a \emph{block} if it is isomorphic to
$\cV(\epsilon, \ell)$ for some $\epsilon \in \kut $, $\ell \in \N_{\ge 2}$.

\begin{theorem}\label{th:infGK} \cite[Theorem 1.2]{AAH}
The Gelfand-Kirillov dimension of the Nichols algebra $\NA(\cV(\epsilon,\ell))$ is finite if and only if $\ell = 2$ and $\epsilon^2 =1$. \qed
\end{theorem}

Here is our main result.

\begin{theorem}\label{th:clasif}
Let $H$ be a pointed Hopf algebra, $G = G(H)$ and $V$  its infinitesimal braiding.
Then the following are equivalent:
\begin{enumerate}
\item\label{it:pointed-finiteGK} $\GK H < \infty$ and $V$ is a block.

\item\label{it:pointed-finiteGK-dim2} $\GK H < \infty$, $\dim V = 2$ and $V$ is not of diagonal type.

\item\label{it:clasif} $G$ is nilpotent-by-finite and there exists
a Jordanian or super Jordanian YD-triple $\D = (g, \chi, \eta)$ and $\lambda\in \ku$, $\lambda = 0$ when $\chi^2 \neq \varepsilon$,
such that $V=\cV_g(\chi,\eta)$ and $H \simeq \Ug(\D, \lambda)$, cf. \S \ref{subsubsect:jordan-lifting} and \ref{subsubsect:super-jordan-lifting}.
\end{enumerate}
\end{theorem}

\pf \eqref{it:pointed-finiteGK} $\Rightarrow$ \eqref{it:pointed-finiteGK-dim2}: by Theorem \ref{th:infGK}, $V \simeq \cV(\epsilon,2)$ with $\epsilon^2 =1$, thus $\dim V = 2$ and $V$ is not of diagonal type.
\eqref{it:pointed-finiteGK-dim2} $\Rightarrow$ \eqref{it:clasif}:
by Gromov's theorem, $G$ is nilpotent-by-finite. By Lemma \ref{lema:indec2-coss}, $V$ is a block, hence
Propositions \ref{prop:lifting-jordan} and \ref{prop:lifting-super-jordan} apply;
these Propositions also provide \eqref{it:pointed-finiteGK} $\Leftarrow$ \eqref{it:clasif}.
\epf

The isomorphism classes of the Hopf algebras $\Ug(\D, \lambda)$ are also determined in Propositions \ref{prop:lifting-jordan} and \ref{prop:lifting-super-jordan}.

The paper is organized as follows. In Section \ref{sec:yd-dim2} we recall the definitions of the Nichols algebras called the Jordan and super Jordan planes. We then
 discuss indecomposable Yetter-Drinfeld modules of dimension 2 over groups. Section \ref{section:pre-nichols} is dedicated to a discussion on the problem of generation in degree one, that is equivalent to the study of post-Nichols algebras with finite $\GK$. We show how to reduce (in general) this problem to the study of pre-Nichols algebras with finite $\GK$ (see the relevant definitions below) and deduce from results in \cite[Section 4]{AAH} that the only post-Nichols algebra of the Jordan, or super Jordan, plane with finite $\GK$ is the Nichols algebra itself. Finally we describe all possible liftings of the Jordan, respectively super Jordan, plane in
Proposition \ref{prop:lifting-jordan}, respectively \ref{prop:lifting-super-jordan}.

\subsection{Notation} We refer to \cite{AS Pointed HA} for unexplained terminology and notation.
If $G$ is a group, then $\widehat{G}$ denotes its group of characters.

\section{Yetter-Drinfeld modules of dimension 2}\label{sec:yd-dim2}

\subsection{The Jordan and super Jordan planes}\label{subsec:planes}
We assume from now on that $\epsilon^2 = 1$.
Keep the notation above and set
$x_{21} = \ad_c x_2\, x_1 = x_2x_1 - \epsilon x_1 x_2$.

The Nichols algebra $\cB(\cV(1,2))$ is a  well-known quadratic algebra,
the so-called Jordan plane, related to the quantum Jordan $SL(2)$; it also appears in the classification of AS-regular graded algebras of global dimension 2 \cite{AS}.

In turn, we call $\cB(\cV(-1,2))$ the \emph{super Jordan plane}.

\begin{prop} \label{pr:1block} \cite[Propositions 3.4 \& 3.5]{AAH} The algebras $\cB(\cV(\epsilon,2))$
	have GK-dim 2 and are presented by generators $x_1$ and $x_2$ with defining relations
	\begin{align}\label{eq:rels B(V(1,2))}
	&x_2x_1-x_1x_2+\frac{1}{2}x_1^2,& && &\text{if } \epsilon = 1;
	\\ \label{eq:rels-B(V(-1,2))-1}
	 	&x_2x_{21}- x_{21}x_2 - x_1x_{21}, & &x_1^2,& &\text{if } \epsilon = -1.
	\end{align}
	Further, $\{ x_1^a x_2^b: a,b\in\N_0\}$, respectively $\{x_1^a x_{21}^bx_2^c: a\in\{0,1\},b,c\in\N_0\}$, is a basis of
	$\cB(\cV(1,2))$, respectively $\cB(\cV(-1,2))$.  \qed
	\end{prop}

\subsection{Indecomposable modules over abelian groups}\label{subsection:indec-2}
Let $\Gamma$ be an abe\-lian group. Let
$g\in \Gamma$,  $\chi\in \widehat\Gamma$ and $\eta: \Gamma \to \ku$ a $(\chi, \chi)$-derivation, i.~e.
\begin{align*}
\eta(ht) &= \chi(h)\eta(t) + \eta(h)\chi(t),& h,t &\in \Gamma.
\end{align*}
Let $\cV_g(\chi, \eta) \in \ydG$ be a vector space of dimension 2,
homogeneous of degree $g$ and with action of $\Gamma$ given in a basis $(x_i)_{i\in\I_2}$ by
\begin{align}\label{equation:basis-triangular-gral}
h\cdot x_1 &= \chi(h) x_1,& h\cdot x_2&=\chi(h) x_2 + \eta(h)x_{1},
\end{align}
for all $h \in\Gamma$.
Then $\cV_g(\chi, \eta)$ is indecomposable in $\ydG \iff \eta\neq 0$.
As a braided vector space, $\cV_g(\chi, \eta)$ is either of diagonal type, when $\eta(g) =0$, or else isomorphic to $\cV(\epsilon,2)$,
$\epsilon = \chi(g)$ (indecomposability as Yetter-Drinfeld module is not the same as indecomposability as braided vector space).

\begin{lemma}\label{lema:indec2-abgroup}
Let $V\in \ydG$, $\dim V =2$. Then either $V$ is of diagonal type or else
$V \simeq \cV_g(\chi, \eta)$ for unique $g$, $\chi$ and $\eta$  with $\eta(g) = 1$.
\end{lemma}

\pf Assume that $V$ is not of diagonal type; then $V$ is indecomposable. Since
$\ku \Gamma$ is cosemisimple, there exists $g \in \Gamma$ such that $V$
is homogeneous of degree $g$. Moreover, $\ku =\overline{\ku}$ implies that
$V$ is not simple. Hence there exist $\chi_1, \chi_2 \in \widehat\Gamma$ such
that
$\soc V \simeq \ku_{g}^{\chi_1}$ and $V/\soc V \simeq \ku_{g}^{\chi_2}$.
Pick $x_1\in \soc V -0$ and $x_2\in V_{g_2} - \soc V$;  then
$h\cdot x_2 = \chi_2(h) x_2 + \eta(h)x_{1}$ for all $h \in\Gamma$, where $\eta$ is a $(\chi_1, \chi_2)$-derivation.
Since  $V$ is not of diagonal type,  $ \chi_1(g) =  \chi_2(g)$ and $\eta(g) \neq 0$. Now
\begin{align*}
\chi_1(h)\eta(g) + \eta(h)\chi_2(g) &=  \eta(hg) = \chi_1(g)\eta(h) + \eta(g)\chi_2(h)& \Rightarrow
\chi_1(h) &= \chi_2(h)
\end{align*}
for all $h\in \Gamma$.
Finally, up to changing $x_1$, we may assume that $\eta(g) = 1$.
\epf

\subsection{Indecomposable modules over Hopf algebras}\label{subsection:indec2-coss}

Let $K$ be a Hopf algebra with bijective antipode.
A \emph{YD-pair} \cite{AAGMV} for $K$  is a pair
$(g, \chi) \in G(K) \times \Alg(K, \ku)$ such that
\begin{align}\label{eq:yd-pair}
\chi(h)\,g  &= \chi(h\_{2}) h\_{1}\, g\, \Ss(h\_{3}),& h&\in K.
\end{align}
If $(g, \chi)$ is a YD-pair, then  the one-dimensional vector space $\ku_g^{\chi}$,
with action and coaction given by $\chi $ and $g$ respectively, is in $\ydk$.
Conversely, any  $V \in \ydk$ with $\dim V = 1$ is like this, for unique $g$ and $\chi$. If $(g, \chi)$ is a YD-pair, then $g\in Z(G(K))$.

\medbreak If $\chi_1, \chi_2\in \Alg(K, \ku)$, then the space of $(\chi_1, \chi_2)$-deri\-vations is
\begin{align*}
\Der_{\chi_1,\chi_2}(K, \ku) &= \{\eta\in K^*: \eta(ht) = \chi_1(h)\eta(t) + \chi_2(t)\eta(h) \, \forall h,t \in K\}.
\end{align*}

A \emph{YD-triple} for $K$ is a collection $(g, \chi, \eta)$ where
 $(g, \chi)$, is a YD-pair for $K$, cf. \eqref{eq:yd-pair},
$\eta \in \Der_{\chi,\chi}(K, \ku)$, $\eta(g) = 1$ and
\begin{align}\label{eq:YD-triple}
\eta(h) g_1 &= \eta(h\_2) h\_1 g_2 \Ss(h\_3), & h&\in K.
\end{align}

If $K = \ku G$ is a group algebra, then we can think of the collection $(g, \chi, \eta)$ as in $G$, $\widehat G$, $\Der_{\chi,\chi}(G, \ku)$.

\medbreak
Let $(g, \chi, \eta)$ be a  YD-triple for $K$. Let $\cV_g(\chi,\eta)$ be a vector space with a basis $(x_i)_{i\in\I_2}$,
where action   and coaction are given by
\begin{align*}
h\cdot x_1 &= \chi(h) x_1,& h\cdot x_2&=\chi(h) x_2 + \eta(h)x_{1},& \delta(x_i) &= g\otimes x_i,
\end{align*}
$h\in K$, $i\in \I_2$. Then $\cV_g(\chi,\eta) \in \ydk$, the compatibility being granted by \eqref{eq:yd-pair}, \eqref{eq:YD-triple}.
Since $\eta(g)\neq 0$, then $\cV_g(\chi, \eta)$ is indecomposable in $\ydk$.

\begin{lemma}\label{lema:indec2-coss} Let $G$ be a group.
Let $V\in \ydg$, $\dim V =2$.  Then either $V$ is of diagonal type as a braided vector space or else
$V \simeq \cV_g(\chi, \eta)$ for a unique YD-triple $(g, \chi, \eta)$.
\end{lemma}

\pf Assume first that $V = V_{g_1} + V_{g_2}$ as $\ku G$-comodule, with $g_1 \neq g_2 \in G$. Now $ g_2\cdot V_{g_2} = V_{g_2}$, hence $ g_2\cdot V_{g_1} = V_{g_1}$ and similarly $ g_1\cdot V_{g_2} = V_{g_2}$. Thus $V$ is of diagonal type, a contradiction.
Thus we may assume that
$V = V_g$ for some $g\in G$, and Lemma \ref{lema:indec2-abgroup} applies with $\Gamma = \langle g \rangle$, so that $V \simeq \cV_g(\widetilde{\chi}, \widetilde{\eta})$
for some  $\widetilde{\chi}\in \widehat\Gamma$ and $\widetilde{\eta}$ a $(\widetilde{\chi}, \widetilde{\chi})$-derivation. Then  there is a basis $(x_i)_{i\in\I_2}$ where $g$ acts by $A = \left(\begin{matrix}
\epsilon & 1 \\ 0 & \epsilon
\end{matrix}\right)$. But $g \in Z(G)$, hence any $h \in G$ acts by a matrix in the centralizer of
$A = \left\{\left(\begin{matrix}
a & b \\ 0 & a
\end{matrix}\right): a \in \ku^{\times}, b \in \ku \right\}$. In other words,
$V \simeq \cV_g(\chi, \eta)$ for a unique YD-triple $(g, \chi, \eta)$.
\epf

Let $\D = (g, \chi, \eta)$ be a  YD-triple and $\epsilon := \chi(g)$.
 If $\epsilon = 1$, respectively $-1$, then  we say that $\D$ is a \emph{Jordanian},
 respectively \emph{super Jordanian}, YD-triple.

\begin{remark}\label{rem:extension HA iso}
Let $f\in \AuH H$, $g^f=f(g)$, $\chi^f=\chi\circ f^{-1}$, $\eta^f=\eta\circ f^{-1}$.
Then $\D^f = (g^f,\chi^f, \eta^f)$ is a YD-triple and $\chi^f(g^f)=\chi(g)$, $\eta^f(g^f)=\eta(g)$.
Thus, if $\D$ is Jordanian, respectively super Jordanian, then so is $\D^f$.
Let $V^f = \cV_{g^f}(\chi^f,\eta^f)$, with basis $x_1'$, $x_2'$. Then $f$ extends to a Hopf algebra isomorphism
$\widetilde{f}:T(V)\# H \to T(V^f)\# H$ such that $f(x_i)=x_i'$. Let
\begin{align*}
\Aut \D := \{f\in \AuH H: \D^f = D\}.
\end{align*}
Then we have a morphism of groups $\Aut \D \to \AuH(T(V)\# H)$.
\end{remark}

\section{Generation in degree one}\label{section:pre-nichols}

\subsection{A block plus a point} We shall need a result from \cite{AAH} on some braided vector spaces of dimension 3.
Let $\epsilon \in \{\pm 1\}$, $q_{12}, q_{21}, q_{22} \in \kut$ and $a \in \ku$.
Let $V$ be the braided vector space with basis $x_i$, $i\in \I_3$, and braiding
\begin{align}\label{eq:braiding-block-point}
(c(x_i \otimes x_j))_{i,j\in \I_3} &=
\begin{pmatrix}
\epsilon x_1 \otimes x_1&  (\epsilon x_2 + x_1) \otimes x_1& q_{12} x_3  \otimes x_1
\\
\epsilon x_1 \otimes x_2 & (\epsilon x_2 + x_1) \otimes x_2& q_{12} x_3  \otimes x_2
\\
q_{21} x_1 \otimes x_3 &  q_{21}(x_2 + a x_1) \otimes x_3& q_{22} x_3  \otimes x_3
\end{pmatrix}.
\end{align}
The ghost is $\ghost = \begin{cases} -2a, &\epsilon = 1, \\
a, &\epsilon = -1.
\end{cases}$ If $\ghost \in \N$, then  the ghost is \emph{discrete}.

\begin{theorem}\label{thm:point-block} \cite[Theorem 4.1]{AAH}
If $\GK \NA (V) < \infty$, then $V$ is as in Table \ref{tab:finiteGK-block-point}. \qed
\end{theorem}

\begin{table}[ht]
	\caption{Nichols algebras of a block and a point with finite $\GK$}\label{tab:finiteGK-block-point}
	\begin{center}
		\begin{tabular}{|c|c|c|c|c|}
			\hline $q_{12}q_{21}$ & $\epsilon$  & $q_{22}$  & $\ghost$   & $\GK$  \\
			\hline
			$1$ & $\pm 1$ & $1$ or $\notin \G_{\infty}$ & $0$ &    3
			\\ \cline{3-3}\cline{5-5}
			& & $\in \G_{\infty} - \{1\}$ &  &    2
			\\ \cline{2-5}
			& $1$  & $1$  & \small{discrete}   &     $\ghost + 3$
			\\\cline{3-5}
			& & $-1$  & \small{discrete}  &     $2$
			\\ \cline{3-5}
			& & $\in \G'_3$  &  1 &    $2$
			\\
			\cline{2-5}
			& $-1$ & $1$  & \small{discrete} &      $\ghost + 3$
			\\\cline{3-5}
			& & $-1$  &  \small{discrete} &     $\ghost + 2$
			\\ \hline
			$-1$ & $-1$  & $-1$  &   1 &     $2$
			\\\hline
		\end{tabular}
	\end{center}
\end{table}

\subsection{Pre-Nichols versus post-Nichols}\label{subsection:pre-post-nichols}

Let  $V \in \ydk$ finite-dimensional.
A \emph{post-Nichols algebra} of $V$ is a coradically  graded connected Hopf algebra $\E = \oplus_{n\ge 0} \E^n$ in $\ydk$ such that $\E^1 \simeq V$ \cite{AAR}. A fundamental step in the classification of pointed Hopf algebras with finite $\GK$ is the following.

\begin{question}\label{question:postnichols-finite-gkd-gral}
	Assume that $K = \ku G$, with $G$ nilpotent-by-finite.
	If $V \in \ydk$ has $\GK \NA(V) < \infty$, then determine all  post-Nichols algebras $\E$ of $V$ such that $\GK \E < \infty$.
\end{question}

A \emph{pre-Nichols algebra} of $V$ is a graded connected Hopf algebra $\cB = \oplus_{n\ge 0} \cB^n$ in $\ydk$
such that $\cB^1 \simeq V$ generates $\cB$ as algebra \cite{Ma}.
If $\E$ is a post-Nichols algebra of $V$, then
there is an inclusion $\cB(V)\hookrightarrow\E$ of graded Hopf algebras in $\ydk$ and  the graded dual $\E^d$ is a pre-Nichols algebra of $V^*$, and vice versa.

\begin{lemma}\label{lemma:prenichols-finite-gkd-postnichols}
	Let $\cB$ be a pre-Nichols algebra of   $V$,
	$\E = \cB^d$ (recall that $\dim V < \infty$).
	Then  $\GK \E \leq \GK \cB$. If $\E$ is finitely generated, then $\GK \E = \GK \cB$.
\end{lemma}

\pf  Let $W$ be a finite-dimensional
vector subspace of $\E$; without loss of generality, we may assume that $W$ is graded.
Let $\E_n = \sum_{0\le j \le n} W^j$. Now
there exists $m \in \N$ such that $W \subseteq \oplus_{0\le j \le m} \E^j$;
hence
\begin{align*}
\log_n\dim \E_n &\leq \log_n\dim \oplus_{0\le j \le mn} \E^j = \log_n\dim \oplus_{0\le j \le mn} \cB^j \\ & \overset{\clubsuit}{=} \log_n\dim (\oplus_{0\le j \le m} \cB^j)^n \overset{\heartsuit}{\Rightarrow} \GK \E \leq \GK \cB.
\end{align*}
Here, $\oplus_{0\le j \le mn} \cB^j  = (\oplus_{0\le j \le m} \cB^j)^n$ because $V$ generates
$\cB$, hence $\clubsuit$; while $\heartsuit$ follows by the independence of the generators in the definition of $\GK$.

Conversely, assume that $W$ is a finite-dimensional graded vector subspace generating $\E$.
We claim that  $\E_n \supseteq \oplus_{0\le j \le n} \E^j$. Indeed, it suffices to show that
$\E_n \supseteq \E^n$. For, take $x \in \E^n$; then $x = \sum w_{j_1} \dots w_{j_k}$ with
$w_{j_h}\in W$, and
\begin{align*}
n =  \deg w_{j_1} +  \dots + \deg w_{j_k} \geq k \Rightarrow x \in \E_n.
\end{align*}
Hence $\log_n\dim \E_n \geq \log_n\dim \oplus_{0\le j \le n} \E^j = \log_n\dim \oplus_{0\le j \le n} \cB^j$, therefore $\GK \E \geq \GK \cB$.
\epf

\begin{remark}\label{remark:prenichols-finite-gkd-postnichols}
	The inequality in Lemma \ref{lemma:prenichols-finite-gkd-postnichols} might be strict:
	if $\cB = \kk[T]$ a polynomial ring with $\car \kk > 0$,  then $\E$ is the divided power algebra that has  $\GK \E = 0 < 1 = \GK \kk[T]$.
\end{remark}

\begin{question}\label{question:prenichols-finite-gkd-gral}
	If $V \in \ydk$ has $\GK \NA(V) < \infty$, then determine all  pre-Nichols algebras $\cB$ of $V$ such that $\GK \cB < \infty$.
\end{question}

To solve Question \ref{question:prenichols-finite-gkd-gral} for $V$ is  a
first approximation to solve Question \ref{question:postnichols-finite-gkd-gral}
for $V^*$, since it is open whether $\GK \E < \infty$ implies $\GK
\E^d  < \infty$ for a post-Nichols algebra $\E$. However, the next particular
case is  useful. Consider
the poset of pre-Nichols algebras $\mathfrak{Pre}(V)=\{T(V)/I: I\in\mathfrak{S}\}$
with ordering given by the surjections. We say that
$V$ is \emph{pre-bounded} if every chain
\begin{align}\label{eq:chain}
\dots < \cB[3] < \cB[2]  < \cB[1]  < \cB[0] = \toba (V),
\end{align}
of pre-Nichols algebras over $V$ with finite $\GK$, is finite.

\begin{lemma}\label{lemma:postnichols-finite-gkd-prenichols}
	Let $K$ be a Hopf algebra, $V \in \ydk$ finite-dimensional and $\E \in \ydk$ a post-Nichols algebra of $V$ with $\GK \E < \infty$. If $V^*$ is pre-bounded, then $\E$ is finitely generated and $\GK \E = \GK \E^d$.
	In particular, if the only pre-Nichols algebra of $V^*$ with finite $\GK$ is $\toba(V^*)$,
	then $\E = \toba(V)$.
\end{lemma}

\pf First we construct a chain $\E[0] = \toba(V) \subsetneq \E[1] \dots \subseteq\E$ of finitely generated post-Nichols algebras of $V$. Suppose we have build $\E[n]$
and  that $\E \supsetneq \E[n]$ (otherwise we are done by Lemma \ref{lemma:prenichols-finite-gkd-postnichols}). Pick $x \in \E - \E[n]$ homogeneous of minimal degree $m$. Let $W$ be the Yetter-Drinfeld submodule of $\E^m$ generated by $x$ and let
$\E[n + 1]$ be the subalgebra of $\E$ generated by $\E[n] + W$. Clearly $\E[n+1]$ is a Yetter-Drinfeld submodule of $\E$, hence $\E[n+1] \underline{\ot} \E[n+1]$ is a  subalgebra of $\E\underline{\ot} \E$.
By minimality of $m$, $\Delta(\E[n+1]) \subseteq \E[n+1] \underline{\ot} \E[n+1]$. That is,
$\E[n+1]$ is a finitely generated post-Nichols algebra of $V$ with $\GK \E[n+1] < \infty$.
Thus we have a chain \eqref{eq:chain} of pre-Nichols algebras with
$\cB[n] = \E[n]^d$, and  $\GK \cB[n] < \infty$ for all $n$ by Lemma \ref{lemma:prenichols-finite-gkd-postnichols}.
By hypothesis, there is $n$ such that $\E[n] = \E$ and we are done.
Finally, if the only pre-Nichols algebra of $V^*$ with finite $\GK$ is $\toba(V^*)$,
then $\E = \toba(V)$, because there is only one chain \eqref{eq:chain} for $V^*$.
\epf

\subsection{Post-Nichols algebras of the Jordan and  super Jordan planes}

\begin{lemma}\label{lem:pre-Nichols Jordan} Assume that $V$ is associated to either a Jordanian or
	a super Jordanian YD-triple $\D = (g, \chi, \eta)$
The only post-Nichols algebra of $V$ in $\ydk$ with finite  $\GK$ is $\cB(V)$.
\end{lemma}

\pf The dual $V^*$ corresponds to the YD-triple $\D' = (g^{-1}, \chi^{-1}, \eta\circ \Ss)$;
by Lemma \ref{lemma:postnichols-finite-gkd-prenichols}, it is enough to solve the analogous problem for pre-Nichols algebras. Let $\cB$ be a pre-Nichols algebra of $V$
such that $\GK  \cB<\infty$. We have canonical projections $T(V)\twoheadrightarrow \cB \twoheadrightarrow \cB(V)$.
We shall prove that the defining relations of $\cB(V)$ hold in $\cB$.

\emph{Jordan case}.
  Suppose that
$y=x_2x_1-x_1x_2+\frac{1}{2}x_1^2\neq 0$ in $\cB$; note that $y$ is  primitive and $y\in \cB_{g^2}$.
Let $W \in \ydk$ be spanned by the linearly independent primitive elements $x_1$, $x_2$ and
$y$. Then $\cB(W)$ is a quotient of $\cB$, so $\GK\cB(W)<\infty$.
Notice that $W$ satisfies \eqref{eq:braiding-block-point} for
$y=x_3$, $\epsilon=1$, $q_{12}=q_{21}=q_{22}=1$, $a=2$. By Theorem \ref{thm:point-block},
$\GK\cB(W)=\infty$, a contradiction. Then $y=0$ and $\cB=\cB(V)$.

\emph{Super Jordan case}.
Suppose first that $x_1^2\neq 0$ in $\cB$; note that $x_1^2$ is primitive and $x_1^2\in \cB_{g^2}$.
Let $W\in \ydk$ be spanned by the linearly independent primitives $x_1$, $x_2$ and
$x_1^2$. Then $\cB(W)$ is a quotient of $\cB$, so $\GK  \cB(W)<\infty$. But
$W$ satisfies \eqref{eq:braiding-block-point} for $\epsilon=-1$, $x_3=x_1^2$, $q_{12}=q_{21}=q_{22}=1$, $a=-2$.
By Theorem \ref{thm:point-block}, $\GK \cB(W)=\infty$, a contradiction,
so $x_1^2=0$ in $\cB$. Let $r=x_2x_{21}- x_{21}x_2 - x_1x_{21}$. In $T(V) \# K$ we have
\begin{align}\label{eq:coproduct rel - super Jordan}
\Delta(r)&=r\otimes 1+g^3\otimes r+ x_1g^2 \otimes x_1^2-2x_1^2g\otimes x_2.
\end{align}
Assume that $r\neq 0$  in $\cB$. By the preceding and \eqref{eq:coproduct rel - super Jordan}, $r$ is  primitive, and
$r\in \cB_{g^3}$. Let $W'$ be the space spanned the linearly independent primitives $x_1$, $x_2$, $r$.
Since $\cB(W')$ is a quotient of $\cB$, $\GK  \cB(W')<\infty$. But $W'$ satisfies \eqref{eq:braiding-block-point} for $\epsilon=-1$,
$x_3=y$, $q_{12}=q_{21}=q_{22}=-1$, $a=-3$, so $\GK \cB(W')=\infty$ by
Theorem \ref{thm:point-block}, a contradiction.
Therefore $\cB=\cB(V)$.
\epf

\section{Liftings}\label{section:pointed}
Let $G$ be a nilpotent-by-finite group. If $V\in \ydg$, then $T(V)$ is a Hopf algebra in $\ydg$ and we denote
$\T(V) = T(V) \# \ku G$.
We compute all liftings of the Jordan and super Jordan planes $\cV(\epsilon, 2)$ over $\ku G$.
We follow partly the strategy from \cite{AAGMV}, as the coradical  is assumed to be finite-dimensional in \emph{loc. cit.};
instead we use \cite[Theorem 8]{Gunther} being in the pointed context.

\begin{remark}\label{rem:g2 not trivial} If $\D = (g, \chi, \eta)$ is a Jordanian or super Jordanian YD-triple,
then $g^2\neq 1$, since $g^2\cdot x_2=x_2+2\epsilon x_1$.
\end{remark}

\subsection{Liftings of Jordan planes}\label{subsubsect:jordan-lifting}

Let $\D = (g, \chi, \eta)$ be a Jordanian YD-triple for $\ku G$ and $V=\cV_g(\chi,\eta)$. Let $\lambda \in \ku$ be such that
\begin{align}
\label{eq:norm-jordanian} \lambda &=  0,&  &\text{ if }\chi^2\neq \varepsilon.
\end{align}
Let $\Ug = \Ug(\D, \lambda)$ be the quotient of $\T(V)$ by the relation
\begin{align}\label{eq:lifting-jordan}
x_2x_1-x_1x_2+\frac{1}{2}x_1^2 &= \lambda(1-g^2).
\end{align}
Clearly, $\Ug$ is a Hopf algebra quotient of $\T(V)$.
We now show that any  lifting of a Jordan plane over $\ku G$ is
$\Ug(\D, \lambda)$ for some $\D$, $\lambda$.

\begin{prop}\label{prop:lifting-jordan}
$\Ug$ is a pointed Hopf algebra, and a cocycle deformation of $\gr \Ug\simeq \NA(V)\# \ku G$;
it has $\GK  \Ug = \GK  \ku G + 2$.

Conversely, let $H$ be a pointed Hopf algebra with finite $\GK$  such that $G(H) \simeq G$ and
 the infinitesimal braiding of $H$ is isomorphic to $\cV(1, 2)$.
Then $H \simeq \Ug(\D, \lambda)$ for some YD-triple
  $\D$ and $\lambda \in \ku$  satisfying \eqref{eq:norm-jordanian}.

Moreover $\Ug(\D, \lambda) \simeq \Ug(\D', \lambda')$ iff there exist a Hopf
algebra automorphism $f$ of $\ku G$
and $c\in \ku^{\times}$ such that $\D'=f(D)$ and $\lambda = c\lambda'$.
\end{prop}

\pf
First, we claim that there exists a $(\Ug,\NA(V)\# \ku G)$-biGalois object $\cA$,
so that $\Ug$ is a cocycle deformation of $\NA(V)\# \ku G$. Let $\T = \T(V)$.

Let $X$ be the subalgebra of $\T$ generated by $t=(x_2x_1-x_1x_2+\frac{1}{2}x_1^2)g^{-2}$, which is a polynomial algebra in $t$.
Set $A=\T$. The algebra map $f:X\to A$ determined by $f(t)=t-\lambda g^{-2}$ is $\T$-colinear.
Note that $x_2x_1-x_1x_2+\frac{1}{2}x_1^2 -\lambda$ is stable by the action of $\ku G$ because of \eqref{eq:norm-jordanian}.
By \cite[Remark 5.6 (b)]{AAGMV},
\begin{align*}
\cA  & := A/\langle f(X^+)\rangle =\T/\langle x_2x_1-x_1x_2+\frac{1}{2}x_1^2 -\lambda \rangle \\
& \simeq \Big( T(V)/\langle x_2x_1-x_1x_2+\frac{1}{2}x_1^2 -\lambda \rangle \Big) \# \ku G.
\end{align*}
We claim that $\cA\neq 0$, which reduces to prove that
$$ \E=\E(\D,\lambda) := T(V)/\langle x_2x_1-x_1x_2+\frac{1}{2}x_1^2 -\lambda \rangle \neq 0. $$
The algebra map $T(V)\rightarrow \ku$, $x_1\mapsto c$, $x_2\mapsto 1$, where $c^2=2\lambda$, applies $x_2x_1-x_1x_2+\frac{1}{2}x_1^2 -\lambda$
to $0$, so it factors through $\E$, and thus $\E$ is non-trivial. Now \cite[Theorem 8]{Gunther} applies and $\cA$ is a $\NA(V)\# \ku G$-Galois object.
Now there exists a unique (up to isomorphism) Hopf algebra $L=L(\cA,\NA(V)\# \ku G)$ such that $\cA$ is a $(L,\NA(V)\# \ku G)$-biGalois object. But
$L\simeq \Ug$ by a computation as in \cite[Cor. 5.12]{AAGMV}, and the claim follows.

Moreover $\id_\T:\T\to \T=A$ is a section which restricted to $\ku G$ is an algebra map. Arguing as in \cite[Proposition 5.8 (b)]{AAGMV}
and applying \cite[Theorem 4.2, Corollary 4.3]{Sch} we conclude that there exists a section $\gamma:\NA(V)\# \ku G\to \cA$ which restricted
to $\ku G$ is an algebra map.
Then \cite[Proposition 4.14 (b), (c)]{AAGMV} applies and $\gr \Ug\simeq \NA(V)\# \ku G$, so its infinitesimal braiding is Jordanian.
Also, $\GK  \Ug = \GK  \gr \Ug  = \GK  \ku G + 2$ by  \cite[Theorem 5.4]{Z}.

\bigbreak
Conversely, let $H$ be a pointed Hopf algebra with finite $\GK$  such that $G(H) \simeq G$ and
the infinitesimal braiding of $H$ is $V \simeq \cV(1, 2)$.
By Lemma \ref{lema:indec2-coss}, there is a Jordanian YD-triple
$\D = (g, \chi, \eta)$  such that $V=\cV_g(\chi,\eta)$.
Then $\gr H\simeq\cB\#\ku G$ for some post-Nichols algebra over $V$ such that $\GK  \cB<\infty$. By Lemma \ref{lem:pre-Nichols Jordan}, $\cB=\cB(V)$.
In particular $H$ is generated by $H_0=\ku G$ and $H_1$ as an algebra.
Moreover $H_1/H_0\simeq V\# \ku G$, so there exists a surjective Hopf algebra map $\pi:\T\twoheadrightarrow H$ which
identifies $\ku G$ and applies $x_i$ to a $(g,1)$-primitive element $a_i\in H_1\setminus H_0$. As $x_2x_1-x_1x_2+\frac{1}{2}x_1^2$ is a
$(g^2,1)$-primitive element, $\pi(x_2x_1-x_1x_2+\frac{1}{2}x_1^2)\in\ku(1-g^2)$ so there exists $\lambda\in\ku$ satisfying \eqref{eq:norm-jordanian}
such that $a_2a_1-a_1a_2+\frac{1}{2}a_1^2=\lambda(1-g^2)$. Then $\pi$ factors through $\Ug(\D,\lambda)$,
and this map $\Ug(\D,\lambda)\rightarrow H$ is an isomorphism since their associated coradically graded Hopf algebras coincide.

It remains to show the last statement. Let $F:\Ug(\D, \lambda)\to\Ug(\D', \lambda')$ be
an isomorphism of Hopf algebras. Then $F|_{\ku G}$ is an isomorphism of Hopf algebras since $\ku G$ is the coradical.
We may assume that $F_{|\ku G}=\id_{\ku G}$ by Remark \ref{rem:extension HA iso}.
Now $g=g'$ since $\dim \mathcal{P}_{g,1}(\Ug(\D, \lambda))=3$, $\dim \mathcal{P}_{h,1}(\Ug(\D, \lambda))=1$
for all $h\in G(\ku G)$, $h\neq g$, and the same for $\Ug(\D', \lambda')$. As $F(x_i)\in \mathcal{P}_{g,1}(\Ug(\D', \lambda'))$,
\begin{align*}
F(x_i)&= a_i x_1'+b_i x_2'+c_i(1-g) & \mbox{for some } a_i,b_i,c_i & \in\ku, i=1,2.
\end{align*}
As $F(hx_ih^{-1})= h F(x_i)h^{-1}$ for all $h\in G(\ku G)$ we deduce that $b_1=c_1=0$, $c_2=0$ if $\chi\neq\varepsilon$, $b_2=a_1$,
$\chi=\chi'$ and $\eta=\eta'$. Also,
$$ 0=F \left(x_2x_1-x_1x_2+\frac{1}{2}x_1^2 -\lambda(1-g) \right)=(a_1^2\lambda'-\lambda)(1-g), $$
so $a_1^2\lambda'=\lambda$. The other implication is direct.
\epf

\begin{remark} Recall that $G$ is nilpotent-by-finite. If $G$ is torsion-free, then $H =\Ug(\D, \lambda)$ is a domain. Indeed, $\ku G$ is a domain by \cite{B, FS, Mo},
hence $\gr H \simeq \NA(V)\# \ku G$ is a domain, and so is $H$.
To see that $\NA(V)\# \ku G$ is a domain, filter by giving degree 0 to $x_1$ and $G$  and degree 1 to $x_2$;
the associated graded algebra is $S(V)\otimes \ku G$, clearly a domain.
\end{remark}

\subsection{Liftings of super Jordan planes}\label{subsubsect:super-jordan-lifting}

Let $\D = (g, \chi, \eta)$ be a super Jordanian YD-triple for $\ku G$
and $V=\cV_g(\chi,\eta)$. Let $\lambda \in \ku$ be such that
\begin{align}
\label{eq:norm-super-jordanian} \lambda &=  0,&  &\text{ if }\chi^2\neq \varepsilon.
\end{align}
Let $\Ug = \Ug(\D, \lambda)$ be the quotient of $\T(V)$ by the relations
\begin{align}\label{eq:lifting-super-jordan}
x_1^2 &= \lambda(1-g^2),  &   x_2x_{21}- x_{21}x_2 - x_1x_{21} +2\lambda x_2+ \lambda x_1g^2&=0.
\end{align}
We will prove that all liftings of super Jordan planes are the algebras $\Ug(\D, \lambda)$. The proof follows the same steps as for Jordan planes.

\begin{prop}\label{prop:lifting-super-jordan}
$\Ug$ is a pointed Hopf algebra, a cocycle deformation of $\gr \Ug\simeq \NA(V)\# \ku G$; $\GK  \Ug = \GK  \ku G + 2$.
Conversely, let $H$ be a pointed Hopf algebra with finite $\GK$  such that $G(H) \simeq G$ and
the infinitesimal braiding of $H$ is $\simeq \cV(-1, 2)$.

Then $H \simeq \Ug(\D, \lambda)$ for some YD-triple
  $\D$ and $\lambda \in \ku$  satisfying \eqref{eq:norm-super-jordanian}.

Moreover $\Ug(\D, \lambda) \simeq \Ug(\D', \lambda')$ iff there exist a Hopf
algebra automorphism $f$ of $\ku G$
and $c\in \ku^{\times}$ such that $\D'=f(D)$ and $\lambda = c\lambda'$.
\end{prop}

\pf
First, we claim that there exists a $(\Ug,\NA(V)\# \ku G)$-biGalois object $\cA$, so $\Ug$ is a cocycle deformation of $\NA(V)\# \ku G$.
We proceed in two steps.

Let $X_1$ be the subalgebra of $\T$ generated by $t_1 = x_1^2g^{-2}$, which is a polynomial algebra in $t_1$.
Set $A=\T$. The algebra map $f:X_1\to A$ determined by $f(t_1) = t_1-\lambda g^{-2}$ is $\T$-colinear.
Note that $x_1^2 -\lambda$ is stable by the action of $\ku G$ because of \eqref{eq:norm-super-jordanian}.
By \cite[Remark 5.6]{AAGMV},
\begin{align*}
\cA_1  & := A/\langle f(X_1^+)\rangle =\T/\langle x_1^2 -\lambda \rangle \simeq \big( T(V)/\langle x_1^2 -\lambda \rangle \big) \# \ku G.
\end{align*}
We claim that $\cA_1\neq 0$, which reduces to prove that $\E_1=T(V)/\langle x_1^2 -\lambda \rangle \neq 0$.
The algebra map $\psi:T(V)\rightarrow \ku$, $x_1\mapsto c$, $x_2\mapsto 0$, where $c^2=\lambda$,
satisfies $\psi(x_1^2)=\lambda$, $\psi(x_2x_{21}- x_{21}x_2 - x_1x_{21}+2\lambda x_2)=0$. It induces an algebra map
$$ \E=\E(\D,\lambda)= T(V)/\langle x_1^2 -\lambda, x_2x_{21}- x_{21}x_2 - x_1x_{21}+2\lambda x_2 \rangle \to \ku. $$
Thus $\E$ is non-trivial, which implies that $\E_1$ is also non-trivial. Therefore \cite[Theorem 8]{Gunther} applies and $\cA_1$ is a $(T(V)/\langle x_1^2 \rangle)\# \ku G$-Galois object.

In the second step, we consider the subalgebra $X_2$ of $(T(V)/\langle x_1^2 \rangle)\# \ku G$
generated by $t_2 = (x_2x_{21}- x_{21}x_2 - x_1x_{21})g^{-3}$,  a
polynomial algebra in $t_2$.
The algebra map $f:X_2 \to \cA$ determined by $$f(t_2)=(x_2x_{21}- x_{21}x_2 - x_1x_{21}+2\lambda x_2) g^{-3}$$ is $\T$-colinear.
By \cite[Remark 5.6]{AAGMV}, $\cA := \cA_1/\langle f(X_2^+)\rangle =\E \# \ku G$, so $\cA$ is non-trivial.
Then \cite[Theorem 8]{Gunther} applies again and $\cA$ is a
$\NA(V)\# \ku G$-Galois object. Now there exists a unique (up to isomorphism) Hopf algebra $L=L(\cA,\NA(V)\# \ku G)$ such that $\cA$ is a
$(L,\NA(V)\# \ku G)$-biGalois object. But $L\simeq \Ug$ by a computation as in \cite[Cor. 5.12]{AAGMV}, and the claim follows.
Moreover $\gr \Ug\simeq \NA(V)\# \ku G$ and $\GK  \Ug = \GK  \ku G+2$ by \cite[Theorem 5.4]{Z}.

\bigbreak
Conversely, let $H$ be a Hopf algebra such that $H_0 \simeq \ku G$, $\GK  H < \infty$ and the infinitesimal braiding of $H$ is $\simeq \cV(-1, 2)$.
By Lemma \ref{lema:indec2-coss}, there is a super Jordanian YD-triple
$\D = (g, \chi, \eta)$  such that $V=\cV_g(\chi,\eta)$.
Then $\gr H\simeq\cB\#\ku G$ for some post-Nichols algebra over $V$ with $\GK  \cB<\infty$; by Lemma \ref{lem:pre-Nichols Jordan},  $\cB=\cB(V)$.
Arguing as in Proposition \ref{prop:lifting-jordan}, there exists a surjective Hopf algebra map $\pi:\T\twoheadrightarrow H$ which induces a
surjective map $\Ug(\D,\lambda)\rightarrow H$ for some $\lambda$ as in \eqref{eq:norm-super-jordanian},
but $\pi$ is an isomorphism since their associated coradically graded Hopf algebras coincide.

It remains to prove the last statement. Let $F:\Ug(\D, \lambda) \to \Ug(\D', \lambda')$ be a Hopf algebra isomorphism.
As in the proof of Proposition \ref{prop:lifting-jordan} we may assume that $F_{|\ku G}=\id_{\ku G}$, and we have that $g=g'$.
As $F(x_i)\in \mathcal{P}_{g,1}(\Ug(\D', \lambda'))$,
\begin{align*}
F(x_i)&= a_1 x_1'+b_i x_2'+c_i(1-g), & a_i,b_i,c_i & \in\ku, i=1,2.
\end{align*}
Notice that $F(gx_ig^{-1})= g F(x_i)g^{-1}$, so $b_1=c_1=c_2=0$, $b_2=a_1$ since $\chi(g)=\chi'(g)=-1$. But as
$F(hx_ih^{-1})= h F(x_i)h^{-1}$ for all $h\in G(\ku G)$ we conclude that $\chi=\chi'$ and $\eta=\eta'$. Also,
$$ 0=F\left(x_1^2 -\lambda(1-g) \right)=(a_1^2\lambda'-\lambda)(1-g), $$
so $a_1^2\lambda'=\lambda$. The other implication is direct.
\epf

\begin{remark} The algebra $H =\Ug(\D, \lambda)$ is never a domain:
$ab =0$, where
\begin{align*}
a &= \sqrt{\lambda}\,(g - 1) + x,& b &= \sqrt{\lambda}\, (g + 1) + x.
\end{align*}
\end{remark}

\end{document}